
\input gtmacros
\input pictex

\input gtmonout
\volumenumber{2}
\volumeyear{1999}
\volumename{Proceedings of the Kirbyfest}
\pagenumbers{455}{472}
\papernumber{22}
\received{16 November 1998}\revised{8 July 1999}
\published{21 November 1999}

\def\dim{{\rm dim}}
\def\min{{\rm min}}
\def\lk{{\rm lk}}
\def\st{{\rm st}}
\def\im{{\rm im}}
\def\inv{^{-1}}
\def\Z{\Bbb Z}

\def\mbox#1{\hbox{#1}}
\newdimen\unitlength
\def\makebox(0,0)[l]#1{\hbox to 0pt{#1\hss}}
\def\symbol#1{\char#1}
\font\tencirc=lcircle10
\def\SetFigFont#1#2#3#4{{\small$#4$}}
\font\thinlinefont=cmr5

\font\spec=cmtex10 scaled 1095 
\def\d{\hbox{\spec \char'017\kern 0.05em}}       
\font\specs=cmtex10 
\def\smalld{\hbox{\specs \char'017\kern 0.05em}} 

\reflist

\key{BRS} {\bf Sandro Buonchristiano}, {\bf Colin Rourke}, 
{\bf Brian Sanderson}, {\it A geometric approach to homology theory,
VII: the geometry of CW complexes}, LMS Lecture Note Series, 18 (1976)
131--149

\key{GM1} {\bf Mark Goresky}, {\bf Robert MacPherson}, {\it Intersection 
homology theory}, Topology, 19 (1980) 135--162
 
\key{GM2} {\bf Mark Goresky}, {\bf Robert MacPherson}, {\it Intersection 
homology theory: II}, Inventiones Mathematicae, 71 (1983) 77--129
 
\key{GM3} {\bf Mark Goresky}, {\bf Robert MacPherson},  {\it Problems and 
bibliography on intersection homology}, from: ``Intersection cohomology'',
Progress in Mathematics, Birkh\"auser, Boston, 60 (1984) 221--233 
 
\key{King} {\bf Henry C King}, {\it Topological invariance of Intersection 
homology without\break sheaves}, Topology and its Applications, 
20 (1985) 149--160

\key{H-W} {\bf P\,J Hilton}, {\bf S Wylie}, {\it An introduction to 
algebraic topology}, CUP, Cambridge (1965)

\key{Hughes} {\bf Bruce Hughes}, {\it Geometric topology of stratified
spaces}, Electronic Research Announcements of the AMS, 2 (1996) 73--81

\key{Frank} {\bf Frank Quinn}, {\it Homotopically stratified sets}, 
J. Amer. Math. Soc. 1 (1988) 441--499

\key{R-S} {\bf Colin Rourke}, {\bf Brian Sanderson}, {\it         
Introduction to piecewise linear topology}, Springer study edition,
Springer--Verlag, Berlin (1982)
 
\key{Larry} {\bf L\,C Siebenmann}, {\it Deformations of homeomorphisms
of stratified sets}, Comm. Math. Helv. 47 (1972) 123--163

\endreflist

\title{Homology stratifications and intersection homology}
\author{Colin Rourke\\Brian Sanderson}

\address{Mathematics Institute, University of Warwick\\Coventry CV4 7AL, UK}
\email{cpr@maths.warwick.ac.uk, bjs@maths.warwick.ac.uk}

\abstract
A homology stratification is a filtered space with local homology
groups constant on strata.  Despite being used by Goresky and
MacPherson [\ref{GM2}] in their proof of topological invariance of
intersection homology, homology stratifications do not appear to have
been studied in any detail and their properties remain obscure.  Here
we use them to present a simplified version of the Goresky--MacPherson
proof valid for PL spaces, and we ask a number of questions.  The
proof uses a new technique, homology general position, which sheds
light on the (open) problem of defining generalised intersection
homology.
\endabstract

\asciiabstract{%
A homology stratification is a filtered space with local homology
groups constant on strata.  Despite being used by Goresky and
MacPherson [Intersection homology theory: II, Inventiones
Mathematicae, 71 (1983) 77-129] in their proof of topological
invariance of intersection homology, homology stratifications do not
appear to have been studied in any detail and their properties remain
obscure.  Here we use them to present a simplified version of the
Goresky--MacPherson proof valid for PL spaces, and we ask a number of
questions.  The proof uses a new technique, homology general position,
which sheds light on the (open) problem of defining generalised
intersection homology.}

\primaryclass{55N33, 57Q25, 57Q65}

\secondaryclass{18G35, 18G60, 54E20, 55N10, 57N80, 57P05}

\keywords{Permutation homology, intersection homology, homology\break 
stratification, homology general position}
\asciikeywords{Permutation homology, intersection homology, homology 
stratification, homology general position}

\makeshorttitle

\quote{\small\it  Rob Kirby has been a great source of encouragement.  
His help in founding the new electronic journal \gt\ has been invaluable.
It is a great pleasure to dedicate this paper to him.}\endquote

\section{Introduction}

Homology stratifications are filtered spaces with local homology
groups constant on strata; they include stratified sets as special
cases.  Despite being used by Goresky and MacPherson [\ref{GM2}] in
their proof of topological invariance of intersection homology, they
do not appear to have been studied in any detail and their properties
remain obscure.  It is the purpose of this paper is to publicise these
neglected but powerful tools.  The main result is that the
intersection homology groups of a PL homology stratification are given
by singular cycles meeting the strata with appropriate dimension
restrictions.  Since any PL space has a natural intrinsic
(topologically invariant) homology stratification, this gives a new
proof of topological invariance for intersection homology, cf
[\ref{King}].  This new proof is in the spirit of the original proof
of Goresky and MacPherson [\ref{GM2}] who used a similar, but more
technical, definition of homology stratification.  It applies only to
PL spaces, but these include all the cases of serious application (eg
algebraic varieties).  In the proof we introduce a new tool: a
homology general position theorem for homology stratifications.  This
theorem sheds light on the (open) problem of defining intersection
bordism and, more generally, generalised intersection homology.

The rest of this paper is arranged as follows.  In section 2 we define
{\sl permutation homology groups}.  These are groups $H^\pi_i(K)$
defined for any principal $n$--complex $K$ and permutation $\pi\in
\Sigma_{n+1}$.  Permutation homology is a convenient device (implicit 
in Goresky and MacPherson [\ref{GM1}]) for studying intersection
homology.  We prove that, for a PL manifold, all permutation homology
groups are the same as ordinary homology groups.  In section 3 we
prove that the groups are PL invariant for {\sl allowable} permutations by
giving an equivalent singular definition (for a stratified set).  This
makes clear the connection with intersection homology.  In section 4
we extend the arguments of section 2 to homology manifolds and in
section 5 we define homology stratifications, extend the arguments of
sections 3 and 4 to homology stratifications and deduce topological
invariance.  In section 6 we discuss the problem of defining
intersection bordism (and more generally, generalised intersection
homology) in the light of the ideas of previous sections.  Finally in
section 7 we ask a number of questions about homology stratifications.

\section{Permutation homology}

We refer to [\ref{R-S}] for details of PL topology.  Throughout the
paper a {\sl complex} will mean a locally finite simplicial complex
and a {\sl PL space} will mean a topological space equipped with a PL
equivalence class of triangulations by complexes.  Let $K$ be a {\sl
principal $n$--complex}, ie, a complex in which each simplex is the
face of an $n$--simplex.  Let $K^{(1)}$ denote the (barycentric) first
derived complex of $K$.  Recall that $K^{(1)}$ is the subdivision of
$K$ with simplexes spanned by barycentres of simplexes of $K$; more
precisely, if we denote the barycentre of a typical simplex $A_i\in K$
by $a_i$ then a typical simplex of $K^{(1)}$ is of the form
$(a_{i_0},a_{i_1},\ldots,a_{i_k})$ where $A_{i_0}<A_{i_1}<\dots
<A_{i_k}$ and $A_i<A_j$ means $A_i$ is a face of $A_j$.

Now let $\pi\in\Sigma_{n+1}$, the symmetic group, ie, $\pi\co
\{0,1,\ldots,n\}\to\{0,1,\ldots,n\}$ is a permutation.  Define
subcomplexes $K^\pi_i$ of $K^{(1)}$, $0\le i\le n$, to
comprise simplexes $(a_{i_0},a_{i_1},\ldots,a_{i_k})$ where $\dim(A_{i_s})
\in \{\pi(0),\ldots,\pi(i)\}$ for $0\le s\le k$.  In other words
$K^\pi_i$ is the full subcomplex of $K^{(1)}$ generated by the
barycentres of simplexes of dimensions $\pi(0)$, $\pi(1)$ \dots
$\pi(i)$.  The definition implies that $K^\pi_i$ is a principal
$i$--complex and that $K^\pi_i\subset K^\pi_{i+1}$ for each $0\le
i<n$.  Here is an alternative description.  $K^\pi_0$ is the
0--complex which comprises the barycentres of the
$\pi(0)$--dimensional simplexes of $K$ and in general we can describe
$K^\pi_i$ inductively as follows.  To obtain $K^\pi_i$ from
$K^\pi_{i-1}$, attach for each simplex $A_s$ of dimension $\pi(i)$ the
cone with vertex $a_s$ and base $K^\pi_{i-1}\cap \lk(a_s,K^{(1)})$.

\rk{Examples} (cf Goresky and MacPherson [\ref{GM1}, pages 145--147])
\items         
\item{(1)}If $\pi=\hbox{id}$ then $K^\pi_i=K^i$ the $i$--skeleton of $K$.    
\item{(2)}If $\pi(k)=n-k$ for $0\leq k\leq n$ then $K_i^\pi=(DK)^i$ the dual 
$i$--skeleton of $K$.
\item{(3)}For $n=2$ the possibilities for a $2$--simplex intersected with  
$K^\pi_0$ and  $K^\pi_1$ are illustrated in figure 1.         
\figure
\beginpicture\small
\setcoordinatesystem units <0.4cm,0.4cm>
\unitlength=1.04987cm
\linethickness=1pt
\setplotsymbol ({\makebox(0,0)[l]{\tencirc\symbol{'160}}})
\setshadesymbol ({\thinlinefont .})
\setlinear
%
%
\linethickness= 0.500pt
\setplotsymbol ({\thinlinefont .})
\plot  0.476 22.860  2.381 25.718 /
\plot  2.381 25.718  3.810 22.860 /
\putrule from  3.810 22.860 to  0.476 22.860
\putrule from  0.476 22.860 to  0.476 22.860
%
%
\linethickness= 0.500pt
\setplotsymbol ({\thinlinefont .})
\plot  4.763 22.860  6.668 25.718 /
\plot  6.668 25.718  8.096 22.860 /
\putrule from  8.096 22.860 to  4.763 22.860
\putrule from  4.763 22.860 to  4.763 22.860
%
%
\linethickness= 0.500pt
\setplotsymbol ({\thinlinefont .})
\plot  9.049 22.860 10.954 25.718 /
\plot 10.954 25.718 12.383 22.860 /
\putrule from 12.383 22.860 to  9.049 22.860
\putrule from  9.049 22.860 to  9.049 22.860
%
%
\linethickness= 0.500pt
\setplotsymbol ({\thinlinefont .})
\plot 13.335 22.860 15.248 25.709 /
\plot 15.248 25.709 16.669 22.860 /
\putrule from 16.669 22.860 to 13.335 22.860
\putrule from 13.335 22.860 to 13.335 22.860
%
%
\linethickness=1pt
\setplotsymbol ({\makebox(0,0)[l]{\tencirc\symbol{'160}}})
\plot  1.429 24.289  2.191 23.717 /
\plot  2.191 23.717  3.092 24.257 /
%
%
\linethickness=1pt
\setplotsymbol ({\makebox(0,0)[l]{\tencirc\symbol{'160}}})
\putrule from  2.191 23.717 to  2.191 22.860
%
%
\linethickness=1pt
\setplotsymbol ({\makebox(0,0)[l]{\tencirc\symbol{'160}}})
\plot  5.715 24.289  6.477 23.812 /
\plot  6.477 23.812  7.334 24.289 /
%
%
\linethickness=1pt
\setplotsymbol ({\makebox(0,0)[l]{\tencirc\symbol{'160}}})
\putrule from  6.477 23.812 to  6.477 22.860
%
%
\linethickness=1pt
\setplotsymbol ({\makebox(0,0)[l]{\tencirc\symbol{'160}}})
\plot  9.049 22.860 10.954 23.812 /
%
%
\linethickness=1pt
\setplotsymbol ({\makebox(0,0)[l]{\tencirc\symbol{'160}}})
\plot 10.937 25.718 10.954 23.812 /
\plot 10.954 23.812 12.383 22.852 /
%
%
\linethickness=1pt
\setplotsymbol ({\makebox(0,0)[l]{\tencirc\symbol{'160}}})
\plot 17.621 22.860 19.526 25.718 /
\plot 19.526 25.718 21.431 22.860 /
\putrule from 21.431 22.860 to 17.604 22.860
\putrule from 17.621 22.860 to 17.621 22.860
%
%
\linethickness=1pt
\setplotsymbol ({\makebox(0,0)[l]{\tencirc\symbol{'160}}})
\plot 22.384 22.860 24.289 25.718 /
\plot 24.289 25.718 26.194 22.860 /
\putrule from 26.194 22.860 to 22.367 22.860
\putrule from 22.384 22.860 to 22.384 22.860
%
%
\linethickness=1pt
\setplotsymbol ({\makebox(0,0)[l]{\tencirc\symbol{'160}}})
\putrule from 15.240 25.718 to 15.240 23.812
\plot 15.240 23.812 16.669 22.860 /
%
%
\linethickness=1pt
\setplotsymbol ({\makebox(0,0)[l]{\tencirc\symbol{'160}}})
\plot 15.240 23.812 13.335 22.860 /
%
%
\put{\SetFigFont{8}{9.6}{rm}{(02)}} <-1pt,0pt> [lB] at  1.429 21.907
%
%
\put{\SetFigFont{8}{9.6}{rm}{(012)}} <-3pt,0pt> [lB] at  5.715 21.907
%
%
\put{\SetFigFont{8}{9.6}{rm}{(12)}} [lB] at 10.001 21.907
%
%
\put{\SetFigFont{8}{9.6}{rm}{(021)}} <-2pt,0pt> [lB] at 14.287 21.907
%
%
\put{\SetFigFont{8}{9.6}{rm}{\rm Id}} <5pt,0pt> [lB] at 18.574 21.907
%
%
\put{\SetFigFont{8}{9.6}{rm}{(10)}} <2pt,0pt> [lB] at 23.336 21.907
\linethickness=0pt
\putrectangle corners at  0.451 25.764 and 26.240 21.863
\endpicture
\endfigure
\item{(4)}For $n=3$ the intersection of a $3$--simplex with           
$K_0^\pi,\ K^\pi_1$ and $K^\pi_2$ is shown in figure 2 for various $\pi$.
\enditems         
\figure
\beginpicture
\setcoordinatesystem units <.4cm,.4cm>
\unitlength=1.04987cm
\linethickness=1pt
\setplotsymbol ({\makebox(0,0)[l]{\tencirc\symbol{'160}}})
\setshadesymbol ({\thinlinefont .})
\setlinear
%
%
\linethickness=2pt
\setplotsymbol ({\makebox(0,0)[l]{\tencirc\symbol{'161}}})
\ellipticalarc axes ratio  0.085:0.085  360 degrees 
	from  3.768 20.934 center at  3.683 20.934
%
%
\linethickness=3pt
\setplotsymbol ({\makebox(0,0)[l]{\tencirc\symbol{'162}}})
\ellipticalarc axes ratio  0.042:0.042  360 degrees 
	from 19.092 25.252 center at 19.050 25.252
%
%
\linethickness=3pt
\setplotsymbol ({\makebox(0,0)[l]{\tencirc\symbol{'162}}})
\ellipticalarc axes ratio  0.042:0.042  360 degrees 
	from 16.256 19.505 center at 16.214 19.505
%
%
\linethickness=3pt
\setplotsymbol ({\makebox(0,0)[l]{\tencirc\symbol{'162}}})
\ellipticalarc axes ratio  0.042:0.042  360 degrees 
	from 19.547 19.050 center at 19.505 19.050
%
%
\linethickness=3pt
\setplotsymbol ({\makebox(0,0)[l]{\tencirc\symbol{'162}}})
\ellipticalarc axes ratio  0.042:0.042  360 degrees 
	from 22.426 20.500 center at 22.384 20.500
%
%
\linethickness=3pt
\setplotsymbol ({\makebox(0,0)[l]{\tencirc\symbol{'162}}})
\ellipticalarc axes ratio  0.042:0.042  360 degrees 
	from 10.604 21.050 center at 10.562 21.050
%
%
\linethickness=3pt
\setplotsymbol ({\makebox(0,0)[l]{\tencirc\symbol{'162}}})
\ellipticalarc axes ratio  0.042:0.042  360 degrees 
	from 11.229 19.505 center at 11.187 19.505
%
%
\linethickness=3pt
\setplotsymbol ({\makebox(0,0)[l]{\tencirc\symbol{'162}}})
\ellipticalarc axes ratio  0.042:0.042  360 degrees 
	from 12.901 21.442 center at 12.859 21.442
%
%
\linethickness= 0.500pt
\setplotsymbol ({\thinlinefont .})
\plot  3.810 25.241  0.953 19.526 /
\plot  0.953 19.526  4.286 19.050 /
\plot  4.286 19.050  7.144 20.479 /
\plot  7.144 20.479  3.810 25.241 /
%
%
\linethickness= 0.500pt
\setplotsymbol ({\thinlinefont .})
\plot 11.430 25.241  8.572 19.526 /
\plot  8.572 19.526 11.906 19.050 /
\plot 11.906 19.050 14.764 20.479 /
\plot 14.764 20.479 11.430 25.241 /
%
%
\linethickness= 0.500pt
\setplotsymbol ({\thinlinefont .})
\plot 19.050 25.241 16.192 19.526 /
\plot 16.192 19.526 19.526 19.050 /
\plot 19.526 19.050 22.384 20.479 /
\plot 22.384 20.479 19.050 25.241 /
%
%
\linethickness= 0.500pt
\setplotsymbol ({\thinlinefont .})
\plot  3.810 25.241  4.286 19.050 /
%
%
\linethickness= 0.500pt
\setplotsymbol ({\thinlinefont .})
\plot 11.430 25.241 11.906 19.050 /
%
%
\linethickness= 0.500pt
\setplotsymbol ({\thinlinefont .})
\plot 19.050 25.241 19.526 19.050 /
%
%
\linethickness= 0.500pt
\setplotsymbol ({\thinlinefont .})
\plot  2.292 22.257  2.959 21.145 /
\plot  2.959 21.145  4.049 21.897 /
\plot  4.049 21.897  5.245 21.505 /
\plot  5.245 21.505  5.556 22.775 /
\plot  5.556 22.775  4.128 22.267 /
%
%
\linethickness= 0.500pt
\setplotsymbol ({\thinlinefont .})
\plot  4.381 19.622  5.520 19.685 /
\plot  5.520 19.685  5.245 21.495 /
%
%
\linethickness= 0.500pt
\setplotsymbol ({\thinlinefont .})
\plot  3.255 21.960  3.969 22.193 /
%
%
\linethickness= 0.500pt
\setplotsymbol ({\thinlinefont .})
\plot  3.277 21.950  2.303 22.257 /
%
%
\linethickness=1pt
\setplotsymbol ({\makebox(0,0)[l]{\tencirc\symbol{'160}}})
\plot 11.843 21.175 12.848 21.461 /
%
%
\linethickness=1pt
\setplotsymbol ({\makebox(0,0)[l]{\tencirc\symbol{'160}}})
\plot 11.161 20.966 11.627 21.114 /
%
%
\linethickness= 0.500pt
\setplotsymbol ({\thinlinefont .})
\plot  8.594 19.526 10.579 21.050 /
\plot 10.579 21.050 11.436 25.193 /
\plot 11.436 25.193 12.833 21.463 /
\plot 12.833 21.463 14.738 20.479 /
%
%
\linethickness= 0.500pt
\setplotsymbol ({\thinlinefont .})
\plot 12.848 21.431 11.913 19.067 /
%
%
\linethickness= 0.500pt
\setplotsymbol ({\thinlinefont .})
\plot 10.596 21.033 11.927 19.082 /
%
%
\linethickness=1pt
\setplotsymbol ({\makebox(0,0)[l]{\tencirc\symbol{'160}}})
\plot 10.611 21.065 11.168 20.955 /
%
%
\linethickness=1pt
\setplotsymbol ({\makebox(0,0)[l]{\tencirc\symbol{'160}}})
\plot 16.216 19.511 18.804 21.082 /
\plot 18.804 21.082 19.503 19.050 /
%
%
\linethickness=1pt
\setplotsymbol ({\makebox(0,0)[l]{\tencirc\symbol{'160}}})
\plot 22.375 20.496 19.598 20.938 /
%
%
\linethickness=1pt
\setplotsymbol ({\makebox(0,0)[l]{\tencirc\symbol{'160}}})
\plot 18.804 21.065 19.217 20.987 /
%
%
\linethickness=1pt
\setplotsymbol ({\makebox(0,0)[l]{\tencirc\symbol{'160}}})
\plot 18.804 21.050 19.042 25.210 /
%
%
\linethickness= 0.500pt
\setplotsymbol ({\thinlinefont .})
\plot  3.810 19.907  7.144 20.479 /
%
%
\linethickness= 0.500pt
\setplotsymbol ({\thinlinefont .})
\putrule from  3.715 19.622 to  4.096 19.622
%
%
\linethickness= 0.500pt
\setplotsymbol ({\thinlinefont .})
\plot  2.953 21.145  2.762 19.241 /
%
%
\linethickness= 0.500pt
\setplotsymbol ({\thinlinefont .})
\plot  0.953 19.526  2.667 19.812 /
%
%
\linethickness= 0.500pt
\setplotsymbol ({\thinlinefont .})
\plot  2.953 19.812  3.524 19.907 /
%
%
\linethickness= 0.500pt
\setplotsymbol ({\thinlinefont .})
\plot  2.762 19.241  3.715 19.622 /
%
%
\linethickness= 0.500pt
\setplotsymbol ({\thinlinefont .})
\plot  3.715 19.622  3.238 19.876 /
%
%
\linethickness= 0.500pt
\setplotsymbol ({\thinlinefont .})
\plot  3.260 19.865  3.238 20.955 /
%
%
\linethickness= 0.500pt
\setplotsymbol ({\thinlinefont .})
\plot  3.291 21.929  3.260 21.463 /
%
%
\linethickness= 0.500pt
\setplotsymbol ({\thinlinefont .})
\plot  3.281 21.294  3.260 21.188 /
%
%
\linethickness=2pt
\setplotsymbol ({\makebox(0,0)[l]{\tencirc\symbol{'161}}})
\plot  3.704 19.622  3.694 20.934 /
%
%
\linethickness=2pt
\setplotsymbol ({\makebox(0,0)[l]{\tencirc\symbol{'161}}})
\plot  3.662 20.955  2.953 21.135 /
%
%
\linethickness=2pt
\setplotsymbol ({\makebox(0,0)[l]{\tencirc\symbol{'161}}})
\plot  3.683 20.955  4.000 21.061 /
\plot  4.000 21.061  4.011 21.071 /
%
%
\linethickness=2pt
\setplotsymbol ({\makebox(0,0)[l]{\tencirc\symbol{'161}}})
\plot  4.244 21.135  5.228 21.505 /
%
%
\linethickness=2pt
\setplotsymbol ({\makebox(0,0)[l]{\tencirc\symbol{'161}}})
\plot  3.694 20.944  3.503 21.400 /
%
%
\linethickness=2pt
\setplotsymbol ({\makebox(0,0)[l]{\tencirc\symbol{'161}}})
\plot  3.418 21.579  3.291 21.960 /
%
%
\linethickness=2pt
\setplotsymbol ({\makebox(0,0)[l]{\tencirc\symbol{'161}}})
\plot  3.694 20.976  3.662 20.892 /
%
%
\linethickness=2pt
\setplotsymbol ({\makebox(0,0)[l]{\tencirc\symbol{'161}}})
\plot 11.165 20.944 11.176 20.320 /
\plot 11.176 20.320 11.187 20.309 /
%
%
\linethickness= 0.500pt
\setplotsymbol ({\thinlinefont .})
\plot  8.583 19.516 10.996 19.844 /
%
%
\linethickness= 0.500pt
\setplotsymbol ({\thinlinefont .})
\plot 11.292 19.854 11.324 19.865 /
%
%
\linethickness= 0.500pt
\setplotsymbol ({\thinlinefont .})
\plot 11.462 19.897 11.726 19.950 /
%
%
\linethickness= 0.500pt
\setplotsymbol ({\thinlinefont .})
\plot 11.938 19.981 12.171 20.024 /
%
%
\linethickness= 0.500pt
\setplotsymbol ({\thinlinefont .})
\plot 12.435 20.045 14.764 20.479 /
%
%
\linethickness=2pt
\setplotsymbol ({\makebox(0,0)[l]{\tencirc\symbol{'161}}})
\plot 11.176 20.034 11.187 19.505 /
%
%
\linethickness=2pt
\setplotsymbol ({\makebox(0,0)[l]{\tencirc\symbol{'161}}})
\plot 11.176 20.944 11.271 21.304 /
%
%
\linethickness= 0.500pt
\setplotsymbol ({\thinlinefont .})
\putrule from  8.583 19.516 to 11.197 19.516
%
%
\linethickness= 0.500pt
\setplotsymbol ({\thinlinefont .})
\plot 11.187 19.526 11.896 19.071 /
%
%
\linethickness= 0.500pt
\setplotsymbol ({\thinlinefont .})
\plot 11.208 19.516 11.494 19.590 /
%
%
\linethickness= 0.500pt
\setplotsymbol ({\thinlinefont .})
\plot 11.652 19.643 11.790 19.674 /
%
%
\linethickness= 0.500pt
\setplotsymbol ({\thinlinefont .})
\plot 11.959 19.696 12.097 19.738 /
%
%
\linethickness= 0.500pt
\setplotsymbol ({\thinlinefont .})
\plot 12.308 19.727 14.764 20.479 /
%
%
\linethickness= 0.500pt
\setplotsymbol ({\thinlinefont .})
\plot  8.594 19.537 10.583 20.807 /
%
%
\linethickness= 0.500pt
\setplotsymbol ({\thinlinefont .})
\plot 10.742 20.902 10.816 20.934 /
%
%
\linethickness= 0.500pt
\setplotsymbol ({\thinlinefont .})
\plot 10.911 21.093 11.282 21.304 /
%
%
\linethickness= 0.500pt
\setplotsymbol ({\thinlinefont .})
\plot 11.271 21.273 11.441 25.231 /
%
%
\linethickness= 0.500pt
\setplotsymbol ({\thinlinefont .})
\plot 11.282 21.294 11.631 21.220 /
%
%
\linethickness= 0.500pt
\setplotsymbol ({\thinlinefont .})
\plot 11.927 21.103 12.562 20.987 /
%
%
\linethickness= 0.500pt
\setplotsymbol ({\thinlinefont .})
\plot 12.785 20.923 14.743 20.489 /
%
%
\linethickness= 0.500pt
\setplotsymbol ({\thinlinefont .})
\plot 16.214 19.484 19.018 19.960 /
%
%
\linethickness= 0.500pt
\setplotsymbol ({\thinlinefont .})
\plot 19.272 20.003 19.367 20.034 /
\putrule from 19.367 20.034 to 19.357 20.034
%
%
\linethickness= 0.500pt
\setplotsymbol ({\thinlinefont .})
\plot 19.569 20.055 22.394 20.489 /
\linethickness=0pt
\putrectangle corners at  0.927 25.356 and 22.490 18.946
\endpicture
\endfigure

\rk{Definition} The {\sl $i^{th}$ permutation homology group},           
$H_i^\pi(K)$,           
of $K$ is the $i^{th}$ homology group of the chain complex:          
$$\dots\longrightarrow H_{i+1}(K^\pi_{i+1},K^\pi_i)          
\buildrel\smalld\over \longrightarrow H_{i}(K^\pi_{i},K^\pi_{i-1})          
\buildrel\smalld\over \longrightarrow H_{i-1}(K^\pi_{i-1},K^\pi_{i-2})
\longrightarrow\dots$$          
where the boundary homomorphisms come from boundaries in the homology
exact sequencies of the appropriate triples.  Cohomology groups
$H^i_\pi(K)$ are defined similarly.  The definition also extends to
any generalised homology theory;  but see the discussion in section 7.

Using a standard diagram chase (and the fact that homology groups
vanish above the dimension of the complex) we have:

\proc{Proposition}\key{im-prop} $H_i^\pi(K)\cong Im(H_i(K^\pi_i)\to         
H_i(K^\pi_{i+1}))$\qed\endproc          
It follows that $H^\pi_i(K)$ can be described as $i$--cycles
in $|K^\pi_i|$ modulo homologies in $|K^\pi_{i+1}|$          
%
and we are at liberty to use singular or simplicial cycles and
homologies.  By releasing the restriction on cycles and boundaries we
get a natural map $\phi\co H^\pi_i(K)\to H_i(K)$.

\proc{Proposition}\key{PL-man}
If $|K|$ is a PL manifold then the natural map $\phi\co
H^\pi_i(K)\to H_i(K)$ is an isomorphism.\endproc

\prf 
The vertices of $K^{(1)}$ not used in the construction of $K^\pi_i$
consist of barycentres of simplexes $A$ with $\dim(A)\not\in\pi[0,i]$
and we denote by $CK^\pi_i$ the full subcomplex (of dimension $n-i-1$)
generated by these unused vertices.  This can also be defined as
follows: write $\bar\pi(k)=n-\pi(k)$ then $CK^\pi_i:=
K^{\bar\pi}_{n-i-1}$.  Note that $|K^\pi_i|\,\cap\,|CK^\pi_i| =\emptyset$
and any simplex of $K^{(1)}$ may be uniquely expressed as a join of a
simplex of $K^\pi_i$ with a simplex of $CK^\pi$. Now an $i$--cycle in
$|K|$ may be pushed off $|CK^\pi|$ by general position and then it can
be pushed down join lines into $|K^\pi_i|$. Similarly homologies can
be pushed off $|CK^\pi_{i+1}|$ into $|K^\pi_{i+1}|$.
\endprf

\section{PL invariance}\key{sec:PL-inv}

Now let $d^\pi_{i,j}$ be $|\pi[0,i]\cap [0,j]|-1$, ie, one less
than the number of integers $\le i$ which have image under $\pi$ which
is $\le j$.

The following facts are readily checked:

\proc{Lemma}\key{d-lem}\items

\item{\rm(1)}The integers $d^\pi_{i,j}$ satisfy $d^\pi_{i,j}\le\min(i,j)$,
$d^\pi_{n,j}=j$, $d^\pi_{i,n}=i$,  $d^\pi_{i,j}-d^\pi_{i-1,j}=0$ or 1,
$d^\pi_{i,j}-d^\pi_{i,j-1}=0$ or 1.

\item{\rm(2)}The integers $d^\pi_{i,j}$ determine the permutation $\pi$.

\item{\rm(3)}$d^\pi_{i,j}$ is the dimension of $K^\pi_i\cap K^j$ where
$K^j$ is the $j$--skeleton of $K$.\qed\enditems\endproc

We now use the integers $d^\pi_{i,j}$ to define {\sl singular permutation
homology} for a filtered space.

Define a {\sl (geometric) $n$--cycle} (often called a {\sl
pseudo-manifold}) to be a compact oriented PL $n$--manifold with
singularity of codimension $\ge2$.  This is the natural picture for a
(glued-up) singular cycle.  A {\sl cycle with boundary} is a compact
oriented PL manifold with boundary and singularity of codimension
$\ge2$, which meets the boundary in codimension $\ge2$.  In other
words if $P$ is a cycle with boundary then $\d P$ is a cycle of one
lower dimension.  By a {\sl (geometric) singular cycle} $(P,f)$ in a
space $X$ we mean a geometric $n$--cycle $P$ and a map $f\co P\to X$.
A {\sl (geometric) singular homology} $(Q,F)$ between singular cycles
$(P,f)$, $(P',f')$ is a cycle $Q$ with boundary isomorphic to $P\cup
-P'$ such that $F|P=f$, $F|P'=f'$.  It is well known that (singular)
homology can be described as geometric singular homology classes of
geometric singular cycles.
There is a similar description for relative singular homology.  A {\sl
relative} singular cycle $(P,f)$ in a pair of spaces $(X,A)$ is a
geometric cycle $P$ with boundary $\d P$ and a map of pairs $f\co (P,\d P)
\to (X,A)$.  A {\sl relative homology} $(Q,F)$ between relative cycles
is a cycle $Q$ with boundary isomorphic to $P\cup -P'\cup Z$, where $Z$
is a cycle with boundary $\d P \cup -\d P'$, and $F$ is a map
of pairs $(Q,Z)\to(X,A)$ such that $F|P=f$, $F|P'=f'$.  We shall
refer to $Z$ as the homology {\sl restricted to the boundary}. 
From now singular cycles and homologies will
all be geometric and we shall omit ``geometric''.

Let $\bar X= \{X_0\subset X_1\subset\ldots\subset X_n\}$ be a filtered
space where $X_j$ has (nominal) dimension $j$.  We refer to
$X_j-X_{j-1}$ as the $j^{th}$ {\sl stratum} of $\bar X$ even though we
are not assuming that $\bar X$ is a stratified set and we often
abbreviate $X_n$ to $X$.  Define the {\sl singular permutation homology
group} $SH^\pi_i(\bar X)$ to be the group generated by singular
$i$--cycles $(P,f)$ in $X$ such that $f\inv(X_j)$ is a PL subset of
dimension $\le d^\pi_{i,j}$ modulo homologies $(W,F)$ such that
$F\inv(X_j)$ is a PL subset of dimension $\le d^\pi_{i+1,j}$.  There
is a similar definition of {\sl relative} singular permutation
homology groups.

\proc{Remark}\key{PL-approx}\rm
If $X$ is a PL space filtered by PL subsets then there is no loss in
assuming that the maps $f$ and $F$ in the definition are PL.  This is
because any map can be approximated by a PL map and it can be checked
that this can be done preserving the (PL) preimages of the closures of
the strata.\fnote{In the standard proof of the simplicial
approximation theorem [\ref{H-W}, pages 37--39], suppose that $f\co
K\to L$ is a map such that $f\inv(L_0)=K_0$ (subcomplexes).  By
subdividing if necessary assume that $L_0$ is a full subcomplex of
$L$.  Suppose that $K$ is sufficiently subdivided for the simplicial
approximation to be defined.  When constructing the simplicial
approximation $g$, choose images of vertices not in $K_0$ to be not in
$L_0$ then $g\inv(L_0)=K_0$.}
\endrk

A permutation $\pi$ is {\sl allowable} if the integers $d^\pi_{i,j}$
satisfy the further condition:
$$d^\pi_{i+1,j}=d^\pi_{i,j}+1\quad\hbox{if}\quad0\le d^\pi_{i,j}<j\eqno{(*)}$$
We shall see that intersection homology groups are precisely the groups
$SH^\pi_i$ for allowable $\pi$.

More generally if $\bar X$ is a filtered space, define $\pi$ to
be $\bar X$--{\sl allowable} if $(*)$ holds for all $j$ such that
$X_j-X_{j-1}\ne\emptyset$. 

It can readily be verified that singular permutation homology has an
excision property for allowable permutations (proved by cutting
cycles and homologies along codimension 1 subsets---allowability is
needed so that the ``constant'' homology is a homology in $SH^\pi$).

Now recall that any PL space $X$ (of dimension $n$) has a natural PL
stratification $\bar X=\{X^0\subset X^1\subset\ldots\subset X^n\}$ where
$X_i-X_{i-1}$ is a PL $i$--manifold.  For any PL stratification
$\bar X$ of $X$, proposition \ref{im-prop} and lemma \ref{d-lem}
provide a natural map $\psi\co H^\pi_i(X)\to SH^\pi_i(\bar X)$.

The following theorem generalises theorem \ref{PL-man} and implies PL
invariance for allowable permutations.

\proc{Theorem}\key{PL-invariance}$\psi\co H^\pi_i(X)\to SH^\pi_i(\bar X)$
is an isomorphism where $\bar X$ is any PL stratification
of $X$ and $\pi$ is $\bar X$--allowable.\endproc

\prf 
To see that $\psi$ is onto we generalise the proof of \ref{PL-man}.
Triangulate $X$ by $K$ say and let $(P,f)$ be a singular $i$--cycle
representing an element of $SH^\pi_i(X)$.  By remark \ref{PL-approx} we
may assume that $f$ is PL; then working inductively over the strata of
$X$ we push $\im(f)$ off $|CK^\pi_i|$ (and hence into $|K^\pi_i|$)
using general position and extending to higher strata using the local
product structure of the stratification.  Notice that the condition
that $\pi$ is $\bar X$--allowable is needed to ensure that the
homologies given by these moves have the correct dimension
restrictions. A similar argument (applied to homologies) shows that
$\psi$ is 1--1.
\endprf

\sh{Connection with intersection homology}

The definition of singular permutation homology is very reminiscent of
the definition of intersection homology.  Indeed we can describe the
connection precisely as follows.  Recall from Goresky and MacPherson
[\ref{GM1}] or King [\King] that a {\sl perversity} is a sequence
$\bar p=\{0=p_0\le p_1\le p_2\le \ldots\le p_n\}$\fnote{Goresky and
MacPherson have the additional condition $p_0=p_1=p_2=0$ and King has
no condition on $p_0$. However if $p_i>i$ then the intersection
condition is vacuous, so we may as well assume $p_0=0$.} where
$p_{i+1}-p_i\le 1.$ Intersection homology (cf [\ref{GM1}, page 138])
is defined exactly like singular permutation homology with
$d^\pi_{i,j}$ replaced by $i+j-n+p_{n-j}$.  However by using
simplicial homology it can be seen that the intersection of an
$i$--cycle with a $j$--dimensional PL subset can always be assumed to
have dimension $\le j$ and a similar remark applies to homologies.
Thus we get exactly the same groups if $d^\pi_{i,j}$ is replaced by
$\min(j\,,\,i+j-n+p_{n-j})$.  We now explain how to find a (unique)
permutation $\pi$ for which $d^\pi_{i,j}$ has this value.

Define a permutation $\pi\in\Sigma_{n+1}$ to be {\sl $V$--shaped} if
$\pi|[0,u]$ is monotone decreasing and $\pi|[u,n]$ is monotone
increasing, where $0\le u\le n$ is the unique number such that
$\pi(u)=0$.  It is easy to see that a $V$--shaped permutation is
uniquely determined by the subset $S_\pi=\pi[0,u-1]
\subset\{1,2,\ldots,n\}$.  We shall see that perversities correspond
to $V$--shaped permutations.  Given a perversity $\bar p$, define
$S=\{j:0< j\le n,\, p_{n-j}=p_{n-j+1}\}$ and consider the $V$--shaped
permutation $\pi$ with $S_\pi=S$.  Then inspecting the graph of $\pi$
it can readily be seen that $d^\pi_{i,j}=\min(j,\,i-q_j)$ where
$q_j=|S_\pi\cap [j+1,n]|$.  But from the definition of $S_\pi$,
$q_j=|k:j<k\le n,\,p_{n-k}=p_{n-k+1}|$, and substituting $c$ for $n-k$
we have $q_j=|c:0\le c<n-j,\,p_{c}=p_{c+1}|= n-j-p_{n-j}$ and hence
$d^\pi_{i,j}=\min(j,\,i+j-n+p_{n-j})$ as required.

It is not hard to see, from graphical considerations, that $V$--shaped
permutations are precisely the same as allowable permutations.  Thus
the singular permutation homology groups for allowable permutations
are precisely the intersection homology groups.  Further it can be
seen that, given an $\bar X$--allowable permutation, there is an
allowable permutation with the same values of $d^\pi_{i,j}$ for all
$j$ such that $X_j-X_{j-1}\ne\emptyset$.  Thus the $\bar X$--allowable
singular permutation groups of $\bar X$ are the intersection homology
groups of $\bar X$.  Thus although permutation homology gives a richer
set of definitions than intersection homology, in the cases where the
groups are PL invariant (which we shall see are the same as the cases
where the groups are topologically invariant) the groups defined are
the intersection homology groups.

In section 5 we will need to consider the permutation $\pi'$ of
$\{0,1,\ldots, n-1\}$ associated to a permutation $\pi$ of
$\{0,1,\ldots, n\}$, defined as follows: remove $0$ from the codomain
of $\pi$ and $\pi^{-1}(0)$ from the domain.  This gives a bijection
between two ordered sets of size $n$.  Identify each with
$\{0,1,\ldots,n-1\}$ by the unique order-preserving bijection.  The
resulting permutation is $\pi'$.  We call $\pi'$ the {\sl reduction}
of $\pi$.  If $\pi$ is allowable then so is $\pi'$ and in terms of
perversities, the operation corresponds to ignoring the final term of
the perversity sequence.  It can be checked that, in terms of the
$d$'s, $\pi'$ is defined by $d^{\pi'}_{i-1,j-1}=d^\pi_{i,j}-1$.

\eject

\section{Homology general position}

Recall that a PL space $M$ is a {\sl homology $n$--manifold} if
$H_i(M,M-x)=0$ for $i<n$ and $H_n(M,M-x)=\Z$ for all $x\in M$ or
equivalently if the link of each point in $M$ is a homology
$(n-1)$--sphere.

The purpose of this section is to generalise proposition \ref{PL-man}
to homology manifolds.

\proc{Proposition}\key{homology-man}
If $M$ is a homology manifold then the natural map\nl $\phi\co
H^\pi_i(M)\to H_i(M)$ is an isomorphism.\endproc

The proof is very similar to the proof of \ref{PL-man}.  However the
key point in the proof (the application of PL general position) does
not work in a homology manifold.  In general it is not possible to
homotope a map of an $i$--dimensional set in a homology manifold $M$
off a codimension $i+1$ subset.  However we only need to move off by
a {\it homology} and this can be done.

\proc{Theorem}{\rm (Homology general position)}\qua  Suppose that $M$
is a homology $n$--manifold and $Y\subset M$ a PL subset of dimension
$y$. Suppose that $(P,f)$ is a singular cycle in $M$ of dimension $q$
where $q+y<n$.  Then there is a singular homology $(Q,F)$ between
$(P,f)$ and $(P',f')$ such that $f'(P')\cap Y=\emptyset$.\key{HGP}

Furthermore the ``move'' can be assumed to be arbitrarily small in the
sense that $F(Q)$ is contained within an arbitrarily small
neighbourhood of $f(P)$.\endproc

There is a version of the theorem which applies to cycles with boundary:

\proclaim{Addendum}
Suppose that $P$ has boundary $\d P$ then there is a relative singular
homology $(Q,F)$ between $(P,f)$ and $(P',f')$ such that $f'(P')\cap
Y=\emptyset$.  Further the moves on both $P$ and $\d P$ can be assumed
to be small, ie, $F(Q)$ is contained within an arbitrarily small
neighbourhood of $f(P)$ and $F(Z)$ is contained within an arbitrarily
small neighbourhood of $f(\d P)$ where $Z$ is the restriction of the
homology to the boundary.\endproc

There is also a relative version of the theorem, which we leave the
reader to prove:  {\sl If $f(\d P)\cap Y=\emptyset$ then we can assume
that the homology fixes the boundary in the sense that $Z\cong\d
P\times I$ and $F|Z$ is $F$ composed with projection on $\d P$.}

\prf
We observe that if, in the small version of the addendum, $\d
P=\emptyset$ then $Z=\emptyset$ and the addendum reduces to the main
theorem.  Thus we just have to prove the addendum.  (By contrast the
non-small version of the addendum is vacuous, since there is always a
relative homology to the empty cycle!)

The proof of the addendum is by induction on $n$ (this is the {\sl
main induction\/}; there will be a subsidiary induction).  Using the
fact that $M$ is a PL space and $Y$ a PL subset, we may cover $M$ by
cones (denoted $C_i$, with bases denoted $B_i$) with the property that
each $C_i$ is contained in a larger cone $C_i^+$ of the form $C_i\cup
B_i\times I$ and such that $Y\cap C_i^+$ is a subcone.  Furthermore we
can assume that each $C_i^+$ has small diameter and that the $C_i^+$
form $n+1$ disjoint subfamilies, ie, two cones in the same family do
not meet.  This implies that any subset of more than $n+1$ of the
$C_i^+$ has empty intersection.  (This is seen as follows.  Choose a
triangulation $K$ of $M$ such that $Y$ is a subcomplex and let
$K^{(2)}$ be the second derived.  Define the $C_i$ to be small
neighbourhoods of the vertex stars $\st(v_i,K^{(2)})$ for vertices
$v_i\in K^{(1)}$.  Define the $C_i^+$ to be slightly larger
neighbourhoods.  Smallness is achieved by taking $K$ to have small
mesh and the subfamilies correspond to the dimension of the simplex of
$K$ of which $v_i$ is the barycentre.)  Since $M$ is a homology
manifold, the cones $C_i$ are in fact homology $n$--balls and their
bases $C_i$ are homology $(n-1)$--spheres.

We shall ``move'' $(P,f)$ by a series of moves each supported by one
of the cones $C_i^+$ and with the property that if $\d P\cap C_i^+$ is
empty before the move, then it still is after the move.  We number the
subfamilies $1,\ldots,n+1$ and we order the moves so that all the
moves corresponding to cones in the subfamily 1 occur first and then
subfamily 2 and so on. Thus if each $C_i^+$ has diameter smaller than
$\epsilon\over{n+1}$ then the set of moves corresponding to subfamily
$i$ is supported by the $\epsilon\over{n+1}$--neighbourhood of $f(P)$
and the whole move is supported by the $\epsilon$--neighbourhood of
$f(P)$ with similar properties for the restriction to the boundary.
The individual moves are defined by a subsidiary inductive process
which we now describe.

By remark \ref{PL-approx} we may assume that $f$ is PL.  By
compactness of $f(P)$ choose a finite subset ${\cal
C}=\{C_1,C_2,\dots,C_t\}$ of cones so that $\bigcup\cal C$ is a
neighbourhood of $Y\cap f(P)$ and with the order compatible with the
order on the subfamilies. Define $Y_j= Y\cap(C_1\cup\ldots\cup C_j)$.

Suppose that we have already moved $(P,f)$ so that $f(P)\cap
Y_j=\emptyset$ and so that $\bigcup\cal C$ is still a neighbourhood of
$Y\cap f(P)$.  We shall explain how to move $(P,f)$ off $Y$ in
$C=C_{j+1}$ by a move supported in $C^+$ so that $\bigcup\cal C$
remains a neighbourhood of $Y\cap f(P)$ and the property
that $f(P)\cap Y_j=\emptyset$ is preserved.  The result is that
$f(P)\cap Y_{j+1}=\emptyset$.  This inductive process starts trivially
and ends with $P\cap Y_t=P\cap Y=\emptyset$ proving the theorem.

For the induction step we have to move $(P,f)$ off $Y$ in $C$.  We
start by applying (genuine) transversality to $B$.  By transversality
we may assume that $f\inv(B)$ is a bicollared subcomplex $R$ of $P$
of dimension $q-1$ which is therefore a cycle (possibly with boundary)
cutting $P$ into two cycles with boundary $P_0$ and $P_1$ where
$P_1=f^{-1}C$.  Note that $\d P$ also splits at $f\inv(B)$ into two
cycles with boundary $S_0$ and $S_1$ with $\d R=\d S_0=\d S_1$ where
$S_1\subset P_1$.\fnote{The transversality theorem being used here is
elementary.  Projecting onto the collar coordinate we have to make a
PL map $g$ say, from $P$ to an interval, transverse to an interior
point.  But we may assume that $g$ is simplicial and, by inspection, a
simplicial map to an interval is tranverse to all points other than
vertices.  So we just compose $g$ with a small movement in the collar
direction so that $B$ does not project to a vertex.}

We now need to consider two cases.

{\bf Case 1 : $S_1\ne\emptyset$}\qua In this case there is a very easy
move which achieves the required result.  Let $P_1^+$ be a small
neighbourhood of $P_1$ in $P$ and $P_0^-$ the corresponding shrunk
copy of $P_0$.  We ``move'' $(P,f)$ to $(P_0^-,f|)$ by excising
$P_1^+$.  More precisely, we regard $(P\times I, f\circ {\rm proj})$
as a relative homology between $(P,f)$ and $(P_0^-,f|)$ by setting $Z$
(the homology restricted to the boundary) equal to $\d P\times I\cup
P_1^+\times \{1\}$.  If we now let $(P_0^-,f|)$ be the new $(P,f)$
the required properties are clear.

{\bf Case 2 : $S_1=\emptyset$}\qua In this case the easy move
described in case 1 would be fallacious, because we have $\d P\cap
C^+$ non-empty after the move whilst it could well be empty before the
move and the restriction to the boundary of the entire process would
not be small.  We now use the fact that $M$ is a homology manifold.
Since $\d R=\d S_1=\emptyset$, $R$ is a cycle and further $(R,f|)$
bounds $(P_1,f|)$ in $C$.  Since $C$ is a homology ball with boundary
$B$ a homology sphere of dimension bigger than $q-1$, there is a cycle
$(P_2,f_2)$ with boundary $(R,f)$ in $B$ and a cycle with boundary
$(Q,F)$ in $C$ with boundary $(P_1\cup P_2, f|\cup f_2)$.  Extend $Q$
by a collar on $P$ to give a homology between $(P,f)$ and $(P_0\cup_R
P_2, f|\cup f_2)$.  This is the first move.  At this point we use the
main induction hypothesis.  By induction we may make a second move of
$(P_2,f_2)$ off $Y$ in $B$ to $(P'_2,f'_2)$ say.  Using collars this
extends to a move of $(P_0\cup_R P_2, f|\cup f_2)$ to $(P',f')$ say
where ${f'}\inv(B)=P'_2$.  It is clear that $f'(P')\cap Y\cap
C=\emptyset$ and it remains to check that $f'(P')\cap Y_j=\emptyset$
and that $\bigcup\cal C$ is still a neighbourhood of $Y\cap f'(P')$.
But before the start of the induction step $f(P)\cap Y_j=\emptyset$
and since these two are compact they start a definite distance apart;
now the two moves which may have affected this were (1) the
application of genuine transversality to $B$ and (2) the (inductive)
move of $P_2$ off $Y$ in $B$, both of which may be assumed to be
arbitrarily small and hence not affect $f(P)\cap Y_j=\emptyset$.
$\bigcup\cal C$ remains a neighbourhood of $Y\cap f'(P')$ for similar
reasons. $f(P)\cap Y$ starts a definite distance from the frontier of
$\bigcup\cal C$ and the same smallness considerations imply that this
property is preserved.  \endprf

\proof{Proof of proposition \ref{homology-man}}
The analogue of the proof of proposition \ref{PL-man} now proceeds
with obvious changes.  Define $CK^\pi_i$ as before. Then by homology
general position we can move an $i$--cycle in $M$ off $|CK^\pi_i|$ by
a homology and hence by pushing down join lines we can move it into
$|K^\pi_i|$.  Similarly a homology can be moved into
$|K^\pi_{i+1}|$.\endprf

\section{Homology stratifications}

Let $x\in X$ a PL space and let $h$ be any (possibly generalised or
permutation) homology theory.  Then for each $y$ close to $x$
there is a natural map $q\co h_*(X,X-x)\to h_*(X,X-y)$.  This is
because $X-x\to X-\st(x)$ is a homeomorphism where $\st(x)$
denotes a small star of $x$ in $X$.  So define $q\co h_*(X,X-x)\cong
h_*(X,X-\st(x))\buildrel j\over\to h_*(X,X-y)$ where $y\in\st(x)$ and
$j$ is induced by inclusion.

Let $h_*^{\rm loc}(X)$ denote the collection $\{h_*(X,X-x):x\in X\}$
of local homology groups of $X$.
Let $Y\subset X$ define $h_*^{\rm loc}(X)$ to be {\sl locally constant} on $Y$
at $x\in Y$ if $q$ is an isomorphism for $y\in Y$ and $y$ close to $x$.

\rk{Comment} This definition is independent of the PL structure on $X$.
If $X'$ denotes $X$ with a different PL structure then we can find a
star $\st(x,X')\subset\st(x,X)$ and then $q$ factors as
$h_*(X,X-x)\cong h_*(X,X-\st(x,X'))\cong h_*(X,X-\st(x,X))\buildrel
j\over\to h_*(X,X-y)$ and it can be seen that $q$ and $q'$ (the analogous
map for $X'$) coincide.

Further the definition makes sense for a wider class of spaces than PL
spaces---essentially any space with locally contractible
neighbourhoods---for example locally cone-like topologically
stratified sets (Siebenmann's CS sets [\ref{Larry}]).

\rk{Definition} 
A filtered PL space $\bar X= \{X_0\subset X_1\subset\ldots\subset
X_n\}$ is an {\sl $h$--stratification} if $h_*^{\rm loc}(X_n)$ is locally
constant on $X_j-X_{j-1}$ for each $j\le n$.  If $h$ is singular permutation
homology $SH^\pi$ then we call it a {\sl $\pi$--stratification}.
\endrk

A locally trivial filtration with strata homology manifolds (eg a
triangulated CS set) is an $h$--stratification for all $h$.  However
note that $h$--stratifications are weaker than any definition of
topological stratification (eg Hughes [\ref{Hughes}], Quinn
[\ref{Frank}]).  For example a homology manifold (with just one
stratum) is an $h$--stratification for all $h$ but, if not a
topological manifold, is not a topological stratification.  There are
several sensible alternative definitions of homology stratifications,
see the discussion in section 7.

Now any principal complex $X$ of dimension $n$ has an {\sl instrinsic
$h$--stratificat\-ion} defined inductively as follows.  Set $X_n=X$
and define $X_{n-1}$ by $x\not\in X_{n-1}$ if $h_*^{\rm loc}(X)$ is
locally constant at $x$.  If $h_*^{\rm loc}$ is locally constant at a
point in the interior of a simplex $\sigma$ then it is locally
constant on the open star of $\sigma$. It follows that $X_{n-1}$ is a
subcomplex of $X$ of dimension $\le n-1$. In general suppose $X_j$ is
defined.  Define the subcomplex $X_{j-1}\subset X_j$ by $x\not\in
X_{j-1}$ if $x$ is in some $j$--simplex in $X_j$ and $h_*^{\rm
loc}(X)$ is locally constant at $x$ on $X_j$.  It can be seen that
$X_{j}$ is a subcomplex of $X$ of dimension $\le j$.

By definition $\bar X= \{X_0\subset X_1\subset\ldots\subset X_n\}$
is an $h$--stratification.  Further the stratification is topologically
invariant since the conditions which define strata are independent of
the PL structure by the comment made above.

\sh{Topological invariance}

Topological invariance of intersection (ie allowable permutation)
homology is proved by combining the arguments of sections 3 and 4.
The key result follows.

\proc{Main theorem}\key{mainth}
Let $\bar X$ be a $\pi$--stratification where $\pi$ is $\bar X$--allowable.
Then the natural map $\psi\co H^\pi_i(X)\to SH^\pi_i(\bar X)$
is an isomorphism.\endproc

Topological invariance follows at once by applying the theorem to the
(topologically invariant) instrinsic $\pi$--stratification.  The proof
is analogous to the proof of \ref{PL-invariance} and
\ref{homology-man} using the following stratified homology general 
position theorem.

\proc{Theorem}{\rm (Stratified homology general position)}\qua Suppose
that $\bar X$ is a $\pi$--stratification where $\pi$ is
$X$--allowable.  Suppose that $(P,f)$ is a singular $p$--cycle in
$SH^\pi_*(X)$ and suppose that $Y\subset X_n$ is a PL subset such that
$\dim(Y\cap X_j)+ d^\pi_{p,j} <j$ for each $0\le j\le n$. Then there
is a singular homology $(Q,F)$ in $SH^\pi_*(X)$ between $(P,f)$ and
$(P',f')$ such that $f'(P')\cap Y=\emptyset$.

Furthermore the ``move'' can be assumed to be arbitrarily small in the
sense that $F(Q)$ is contained within an arbitrarily small
neighbourhood of $f(P)$.\endproc\key{s-hgp}

The theorem has a version for cycles with boundary analogous to the
addendum to theorem \ref{HGP}:

\proclaim{Addendum}
Suppose that $\bar X$ and $Y$ are as in the main theorem and $(P,f)$
is a singular $p$--cycle with boundary in $X$ which satisfies the
dimension restrictions for a cycle in $SH^\pi_*(X)$.  Then there is a
relative singular homology $(Q,F)$ which satisfies the dimension
restrictions for a homology in $SH^\pi_*(X)$ between $(P,f)$ and
$(P',f')$ such that $f'(P')\cap Y=\emptyset$.  Further the moves on
both $P$ and $\d P$ can be assumed to be small, ie, $F(Q)$ is
contained within an arbitrarily small neighbourhood of $f(P)$ and
$F(Z)$ is contained within an arbitrarily small neighbourhood of $f(\d
P)$ where $Z$ is the restriction of the homology to the
boundary.\endproc

There is also an analogous relative version of the theorem which we
leave the reader to state and prove. 
\prf
The theorem is very similar to the proof of theorem \ref{HGP} with $M$
replaced by $X$ and we shall sketch the proof paying careful attention
only to the places where there is a substantive difference.  We merely
have to prove the addendum and we use induction on $n$.  As before we
may cover $X$ by small cones $C_i\subset C_i^+$ (with the base of
$C_i$ denoted $B_i$) which form $n+1$ disjoint subfamilies and such
that $Y$ meets each in a subcone and such that the local filtration
follows the cone structure.  (In this proof the cones are not homology
balls and the bases are not homology spheres.)

It can be checked that the induced filtration on $B_i$ is a
$\pi'$--stratification; essentially this is because the local homology
of $C_i$ at $B_i$ is the suspension of the local homology of $B_i$.
In the following ``cycle'' means singular cycle in $\pi$ or
$\pi'$--homology as appropriate.

We define a finite subset ${\cal C}=\{C_1,C_2,\dots,C_t\}$ such that
$\bigcup\cal C$ is a neighbourhood of $Y\cap f(P)$ as before and we
set up a subsidiary induction with exactly the same properties.  The
induction proceeds with no change at all for case 1.  For case 2,
which was the first place that properties of $M$ were used, there are
now two subcases to consider.  Let $c$ be the conepoint of $C$ and let
$T$ (a subcone) be the intersection of the stratum of $\bar X$
containing $c$ with $C$.

{\bf Case 2.1}\qua $f(P)\supset T$\qua In this case, by the dimension
hypotheses $Y$ misses $T$ and hence, since $Y$ is a subcone of $C^+$ we
have $Y\cap C^+=\emptyset$, and there is nothing to do.

{\bf Case 2.2}\qua {\sl There is a point $x\in T, x\not\in P$.}\qua In this
case, denote $C-B$ by $C'$.  Now $SH^\pi_*(X,X-x)\cong
SH^\pi_*(X,X-C')$ by the definition of $\pi$--stratification and hence
using excision $SH^\pi_*(C,C-x)\cong SH^\pi_*(C,B)$.  But $(P_1,f|)$
represents the zero class in the former group and hence in the latter.
Thus there is a homology $(Q,F)$ say in $SH^\pi_*$ of $(P_1,f)$ rel
boundary to a class $(P_2,f_2)$ say with $f_2(P_2)\subset B$.  The
proof now terminates exactly as in the previous proof.  We use $(Q,F)$
to move $(P,f)$ to $(P_0\cup_R P_2, f|\cup f_2)$ (the first move) and
then we apply induction to move $(P_2,f_2)$ off $Y$ in $B$ extending
by collars as before to produce $(P',f')$ (the second move).  The
required properties are checked as before. \endprf

\proof{Proof of the main theorem}
The analogue of previous similar proofs now proceeds with obvious
changes.  Triangulate $X$ by $K$ and define $CK^\pi_i$ as before.
Then by stratified homology general position we can move an $i$--cycle
in $SH^\pi(\bar X)$ off $|CK^\pi_i|$ by a homology in $SH^\pi(\bar X)$
and hence by pushing down join lines we can move it into $|K^\pi_i|$.
Similarly a homology can be moved into $|K^\pi_{i+1}|$.\endprf

\section{Intersection bordism}  

We have given three equivalent definitions of permutation homology and
we shall see shortly that there is a hidden fourth definition.  All
four generalise to give definitions of intersection bordism (and more
generally of generalised intersection homology).  Only two are the
same for intersection bordism.  We shall see that these two are
topologically invariant.  

The three equivalent definitions of the $i^{th}$ permutation homology
group were:

\items\item{(1)}The homology of the chain complex:
$$\dots\longrightarrow H_{i+1}(K^\pi_{i+1},K^\pi_i)          
\buildrel\smalld\over \longrightarrow H_{i}(K^\pi_{i},K^\pi_{i-1})          
\buildrel\smalld\over \longrightarrow H_{i-1}(K^\pi_{i-1},K^\pi_{i-2})
\longrightarrow\dots$$          

\item{(2)}Cycles in $K^\pi_{i}$ modulo homologies in $K^\pi_{i+1}$.

\item{(3)}Singular permutation homology of a stratified set, ie, 
singular $i$--cycles meeting strata of dimension $j$ in dimension $\le
d^\pi_{i,j}$ modulo homologies meeting strata of dimension $j$ in
dimension $\le d^\pi_{i+1,j}$.\enditems

The fourth equivalent definition follows from definition (2) using the
property that $K^\pi_{i}$ meets $K^j$ in dimension $\le d^\pi_{i,j}$,
see lemma \ref{d-lem}:

\items\item{(4)}Singular $i$--cycles in $K^\pi_{i}$ which meet $K^j$ in 
dimension $\le d^\pi_{i,j}$ modulo homologies in $K^\pi_{i+1}$ which
meet $K^j$ in dimension $\le d^\pi_{i+1,j}$.\enditems

Now let $h$ denote smooth bordism then we can define permutation
bordism theory (denoted $h^\pi$) in direct analogy to permutation
homology in any of the four ways listed above.

There are natural maps between the four definitions of $h^\pi$ as
follows $(3)\leftarrow(4)\to(2)\to(1)$.  We shall see shortly that
$(4)\to(3)$ is an isomorphism.  There is no reason to expect that
either of $(4)\to(2)\to(1)$ are isomorphisms.  To prove $(2)\to(1)$ is
an isomorphism for homology the fact that homology groups vanish above
the dimension of the complex is used; this is false for bordism.  To
prove that $(4)\to(2)$ is an isomorphism another fact special to
homology is used, namely that a cycle can be assumed to be simplicial
and hence a subcomplex.  Again this is in general false for bordism.
In favour of the two equivalent definitions (3) and (4) we have the
following result.

\proc{Theorem} Definitions $(3)$ and $(4)$ are equivalent for bordism
and define a topological invariant of $X$.\endproc

\proof{Sketch of proof}Stratified homology general position (theorem
\ref{s-hgp}) 
can be extended in two ways (1) replace $\pi$--stratifications by
$h^\pi$--stratifications and $SH^\pi$ by $Sh^\pi$ (ie definition (3)
above) and (2) delete the condition $\dim(Y\cap X_j)+ d^\pi_{p,j} <j$
and alter the conclusion to get $\dim(f'(P')\cap Y\cap X_j)\le
\dim(Y\cap X_j)+ d^\pi_{p,j}-j$.  The proof is the same with obvious
changes.  This implies that a cycle in $Sh^\pi$ can be assumed to meet
$K^j$ in the appropriate dimension by applying the theorem with $Y=K^j$
and then the usual argument (make disjoint from $CK^\pi$ and push into
$K^\pi$) yields a cycle in definition (4).  A similar argument applies
to a homology and this proves that definitions (3) and (4) coincide.

Topological invariance follows by applying this to the instrinsic
$h^\pi$--stratification.\endprf

\rk{Remarks}{\bf1})\qua
Definition (4) is briefly considered by Goresky and MacPherson in
[\ref{GM3}, problem 1].  They do not state topological invariance but
they point out that the definition is unlikely to yield any form of
Poincar\'e duality.  In defence of the definition we would observe
that ordinary bordism has no Poincar\'e duality for manifolds (there
is a duality between bordism and cobordism but none between bordism
groups of complementary dimension).  Thus there is no reason to expect
a definition which generalises bordism of a manifold (intersection
homology generalises ordinary homology of a manifold) to satisfy
Poincar\'e duality.

{\bf 2})\qua Let $h$ be any connected generalised homology theory.
Using the main result of [\ref{BRS}] we can regard $h$ as a
generalised bordism theory (given by bordism classes of maps of
suitable manifolds-with-singularity) and hence we can define
permutation $h$--theory in analogy with permutation bordism as above.
The analogue of the theorem is proved in exactly the same way.
However it must be noted that this definition is dependent on the
particular choice of representation for the theory as bordism with
singularities (which in turn depends on a particular choice of CW
structure for the spectrum).  Thus this construction does not
define $h^\pi$ unambiguously.

\section{Questions about homology stratifications}

The following questions are asked in the spirit of a conference
problem session.  We have no clear idea how hard they are and indeed
some may have simple answers which we failed to notice whilst
writing them. 

The simplest definition of homology stratification is given by using
ordinary (integral) homology.  Call such a stratification an {\sl
$H$--stratification}.  Since, by the stable Whitehead theorem, a
homology equivalence induces isomorphisms of all generalised homology
groups, an $H$--stratification is an $h$--stratification for any
generalised homology $h$.  However this is not clear if $h$ is
intersection (ie allowable permutation) homology.

\proclaim{Question 1}Is an $H$--stratification a $\pi$--stratification 
for allowable $\pi$?  In other words, if the local homology groups are
constant on strata, is the same true for local intersection homology
groups?
\endproc

Question 1 is connected to the problem of characterising maps which
induce isomorphisms of intersection homology groups in terms of
ordinary homology.  Here is a related question.  We say that a map
$f\co X\to Y$ of filtered spaces (of dimensions $n$, $m$ respectively)
{\sl respects} the filtration if $f\inv(Y_{m-k})\subset X_{n-k}$ for
each $k$.  A map which respects the filtration induces a homomorphism
$SH^\pi(X)\to SH^\sigma(Y)$, where $\pi$ is a (repeated) reduction
of $\sigma$ or vice versa, (cf King [\King; page 152]).

\proclaim{Question 2}Suppose that $i\co\bar X\subset \bar Y$ is
an inclusion of filtered spaces which respects the filtration and
induces isomorphisms of all ordinary homology groups for all strata
and closures of strata.  Does it follow that $i$ induces isomorphism
of intersection homology groups?\endproc

Question 1 is also related to the problem of functoriality of
intersection homology [\ref{GM3}, problem 4].  Our main theorem gives
an intrinsic definition of intersection homology namely singular
permutation homology of the intrinsic $\pi$--stratification where
$\pi$ is the appropriate allowable permutation.  By the remarks above
question 2, a map which respects the intrinsic $\pi$--stratification
induces a homomorphism $SH^\pi(X)\to SH^\sigma(Y)$.  This is a
somewhat circular characterisation of maps inducing homomorphisms of
intersection homology, since they are characterised in terms of
intersection homology; it is almost as circular as the characterisation
given in [\ref{GM3}, bottom of page 223].  If question 1 has a
positive answer, then the characterisation becomes rather less
circular: maps which respect the intrinsic $H$--stratification induce
homomorphisms of intersection homology.

\proclaim{Question 3}Is there a good geometric characterisation of maps 
which respect the intrinsic $H$--stratification?  For example is it
sensible to ask for a characterisation in terms of properties of point
inverses?
\endproc

We have remarked that a locally trivial filtration with strata
homology manifolds is an $h$--stratification for all $h$.  The
converse is easily seen to be false: glue three homology balls along a
genuine ball in the boundary; the result is a homology stratification
with the interior of the common boundary ball in one stratum, but is
not necessarily locally trivial along that stratum.  Indeed it is not
clear that the strata of an $H$--stratification must be homology
manifolds.

\proclaim{Question 4}Are the strata of an $H$--stratification homology 
manifolds?  Is the same true of a $\pi$--stratification for allowable
$\pi$?\endproc

We now turn to other (stronger) definitions of homology
stratification.  These all have the property that the strata are
obviously homology manifolds.  Goresky and MacPherson use a somewhat
different definition of $h$--stratification.  Their ``canonical''
$\bar p$--filtration [\ref{GM2}, bottom of page 107] is defined
exactly like our instrinsic $h$--stratification except that instead of
our condition that $h_*^{\rm loc}(X)$ is locally constant on
$X_j-X_{j-1}$ for each $j$ they have two conditions: $h_*^{\rm
loc}(X_j)$ {\sl and} $h_*^{\rm loc}(X-X_j)$ are both locally constant
on $X_j-X_{j-1}$ where $h$ is intersection homology (the latter makes
sense: they are using homology with infinite chains, the second local
homology group is the same as $h_{*-1}(\lk(x,X)-\lk(x,X_j))$).  The
two conditions imply that $h_*^{\rm loc}(X)$ is locally constant.  For ordinary
homology if $h_*^{\rm loc}(X)$ and $h_*^{\rm loc}(X_j)$ are both locally constant then so
is $h_*^{\rm loc}(X-X_j)$.  For intersection homology this is not clear.

\rk{Definitions}A {\sl strong} $h$--stratification is one where  
$h_*^{\rm loc}(X)$ and $h_*^{\rm loc}(X_j)$ are both locally constant
on $X_j-X_{j-1}$ for each $j$.  A {\sl GM--strong} $h$--stratification
is one where $h_*^{\rm loc}(X-X_j)$ and $h_*^{\rm loc}(X_j)$ are both
locally constant on $X_j-X_{j-1}$ for each $j$ (this only makes sense
for geometric theories for which the analogue of infinite chains is
defined).  A {\sl very strong} $h$--stratification is one where
$h_*^{\rm loc}(X_k)$ is locally constant on $X_j-X_{j-1}$ for each $k\ge
j$.\endrk

\proclaim{Question 5}What are the relationships between the definitions?
Are the concepts of strong and GM--strong stratifications distinct?
Are there examples of strong stratifications which are not very strong?
Or indeed examples of stratifications which are not strong?\endproc

\references

\Addresses
\recd

\bye